\documentclass[12pt,a4paper]{article}

\usepackage[english,russian]{babel}
\usepackage{indentfirst}
\usepackage{misccorr}
\usepackage{graphicx}
\usepackage{amsmath}
\usepackage{amssymb}
\usepackage{amsthm}
\usepackage{hyperref}
\usepackage{color}

\newif\ifhide
\hidefalse

\newlength{\wdth}

\newtheorem{rTheorem}{Теорема}

\newtheorem{rRemark}{Замечание}

\newtheorem{rEx}{Пример}

\usepackage{hyperref}
\usepackage{graphicx} 
\usepackage{setspace}
\usepackage{graphicx}
\usepackage{tikz}
\usepackage{bbm}
\newcommand{\RNumb}[1]{\uppercase\expandafter{\romannumeral #1\relax}}

\newif\ifru
\rutrue 

\newif\ifen
\enfalse

\ifru
\title{Пример оптимального эргодического управления в среднем в диффузионной модели страховой компании\footnote{Для обоих авторов данная работа поддержана Фондом развития теоретической физики и математики ``БАЗИС''}
}

\author{А.Ю. Веретенников\footnote{Институт проблем передачи информации им. А.А. Харкевича, Российский университет дружбы народов им. Патриса Лумумбы, email: alexander.veretennikov2011@yandex.ru},\;
\; 
Е.Ю. Ященко\footnote{Московский государственный университет им. М.В. Ломоносова, email: elizaveta.iashchenko@math.msu.ru}}

\fi

\begin{document}

\maketitle

\begin{abstract}
Рассматривается эргодический аналог известной диффузионной модели управления рисками и распределения дивидендов для финансовой компании.  В данном несложном примере любопытно как оптимальные управления, -- которые существуют и которых оказывается бесконечное множество, -- сочетаются с эргодическим уравнением Беллмана. 

~

Ключевые слова: диффузионная модель, эргодическое управление в среднем, 
пример

~

\end{abstract}

\section{Введение, постановка задачи}
Рассматривается эргодический аналог модели оптимального стохастического управления из статьи \cite{Taksar2000}. Такая задача может представлять интерес, как пример применения методов стохастического анализа для формализации проблемы длительного устойчивого функционирования страховой компании в условиях неопределённости и риска. 

Опуская детали обоснования диффузионной аппроксимации, которые можно найти, например, в ссылках к статье \cite{Taksar2000}, рассмотрим сразу идеальную предельную модель движения капитала или резерва страховой компании, принятую в работе \cite{Taksar2000} и других:
\begin{equation}\label{sde1}
dX_t=(\mu-a(X_t)) dt+dW_t, \quad X_0=x >0, \quad 0 \leq a \leq M, \quad M > \mu.
\end{equation} 
Здесь постоянная $\mu>0$ имеет смысл скорости роста страховой премии, получаемой компанией, $a(X_t)$ -- назначаемый управлением фирмы  ``отбор'' части премии, идущий на орграсходы, зарплаты, и т.д. (При этом эта часть может быть, в том числе, и больше $\mu$).  Функция $a(\cdot)$ называется стратегией или марковской стратегией и может быть выбрана компанией с целью максимизировать некоторый функционал, см. ниже. Считаем, что дело происходит в некоем идеальном мире, где компания никогда не разоряется, даже если ее резерв становится отрицательным; в то же время, конечно, любая ``разумная'' стратегия должна стараться ухода на отрицательную полуось, по возможности, не допускать, в то же время, обеспечивая какой-то желаемый, постоянный в среднем доход. Поэтому наложим запрет или правило вида 
\begin{equation}\label{a0}
a(x)\equiv 0, \quad \forall \, x \le  0,
\end{equation} 
означающий, что никакие выплаты ``себе'' при отрицательном резерве невозможны. Также, считаем, что действует ограничение: имеется константа $M>\mu$ такая, что 
\begin{equation}\label{Mmu}
0\le a(x) \le M, \quad \forall \, x \in \mathbb R.
\end{equation} 

При применении стратегии $a(\cdot)$, через $\mathsf E^a$ обозначается математическое ожидание, соответствующее вероятностной мере при данной стратегии. Целью оптимизации является функционал\footnote{В принципе, левая и правая части данной формулы могут зависеть от начальной точки $x$. Однако, для эргодических стратегий это обычно константы.} 
\begin{equation}\label{sup}
\liminf_{T \to \infty} \frac{1}{T} \mathsf{E}_x^{a(\cdot)} \int_0^T a(X_t)\, dt \to \sup_{a(\cdot)} =: r.
\end{equation} 
Второй, дополнительной целью, как всегда в задачах оптимизации, является
отыскание стратегии, на которой этот супремум достигается, или, возможно, почти достигается. 

Интуитивно, в данной модели, конечно, существуют эргодические и неэргодические стратегии. Например, тривиальной неэргодической стратегией является $a(x)\equiv 0, \, \forall x \in \mathbb R$. Немного более разумным кажется применять $a(x)= (\mu - \varepsilon)\, 1(x\ge x_0)$ при некотором $x_0\ge 0$. Этой стратегии соответствует транзиентное поведение с уходом на $+\infty$ и ``выигрышем'' $\mu - \varepsilon$. Однако, более-менее ясно, что максимума на неэргодических стратегиях не достичь. Действительно, поскольку действует запрет (\ref{a0}), то неэргодическое, то есть, транзиентное поведение всегда означает уход на $+\infty$. Опять же, интуитивно, на нестрогом уровне, представляется, что в этом случае компания могла бы получить больше, если бы ``в среднем'' отбирала хотя бы немного больше. 

Тривиальной эргодической стратегией является, в частности, правило 
$$
a(x) = M\, 1(x>x_0), \qquad \text{при любом $x_0\ge 0$}, 
$$
причем для такой стратегии существует (единственная) инвариантная мера процесса $\mu^a$, и маргинальные распределения процесса $\mu_t$ сходятся к ней экспоненциально быстро по метрике полной вариации (см., например, \cite{Ver2000}). Можно отметить, что при $M=2\mu$ такая стратегия обеспечит достижение $r=\mu$, в силу очевидной симметрии относительно $x_0$, причем независимо от выбора $x_0\ge 0$.

Еще одно интуитивное соображение говорит о том, что при постоянном $\mu$ (напомним, что $\mu>0$) вряд ли возможно ``в среднем'' за большое время получить больше, чем $\mu$. Действительно, в таком случае резерв непременно должен был бы становиться отрицательным и уходить на $-\infty$; однако, при отрицательных $x$ действует запрет (\ref{a0}), так что такое поведение невозможно. Поэтому естественной гипотезой кажется, что ответ в этой задаче при любом $M>\mu$ должен быть такой: 
\begin{equation}\label{rmu}
r =\mu.
\end{equation} 
Это, действительно, верно, см ниже теорему 1. Ответ на второй вопрос, -- на каких стратегиях этого супремума можно достичь, -- оказывается, на первый взгляд, немного неожиданным. Таких (эргодических) стратегий бесконечно много, причем описание их может даже не опираться на ``обычное'' в подобных задачах эргодическое уравнение Беллмана (см. \cite{Borkar, Ver}, et al.), хотя и никакого  противоречия с уравнением Беллмана тут нет. Выясняется что, по существу, на любой эргодической стратегии, для которой выполнено некоторое неограничительное техническое условие, достигается искомый максимум $r$. 

В следующем разделе 3 для основного результата приводятся две разных версии доказательства, обе без использования уравнения Беллмана. В разделе 4 это уравнение также выписано, решено, и показано, что на его основе результат оказывается, разумеется, таким же.

\section{Основной результат}

\begin{rTheorem}
Для данном модели выигрыш $r$ равен
\begin{equation}
r=\mu.
\end{equation}
Он достигается, во всяком случае, на любой эргодической стратегии $a(\cdot)$, при которой выполнено условие 
\begin{equation}\label{Xt0}
\lim_{t\to\infty}\frac{\mathsf E^a X_t}{t} = 0.
\end{equation}
\end{rTheorem}
Можно отметить, что на эргодических стратегиях условие (\ref{Xt0}), как правило, выполняется; более того, числитель $\mathsf E^a X_t$ в такой ситуации обычно оказывается ограниченным. Вопрос о том, существуют ли неэргодические стратегии, на которых достигается верхняя грань (\ref{sup}), не обсуждается: кажется, что, в принципе, исключить этого нельзя.

\begin{proof}
{\bf 1 способ.}
При любой эргодической стратегии $a(\cdot)$ в силу уравнения (\ref{sde1}) имеем,
\begin{align*}
\mathsf E^a_x X_T - x = \mathsf E^a \int_0^T (\mu - a(X_s))ds. 
\end{align*}
Разделив обе части на $t$, а затем  устремляя $t$ к бесконечности, получаем в силу условия (\ref{Xt0}), что существует предел
\begin{align*}
\lim_{T\to\infty}T^{-1}\mathsf E^a \int_0^T (\mu - a(X_s))ds = 0. 
\end{align*}
Стало быть, 
$$
\lim_{T\to\infty}T^{-1}\mathsf E^a \int_0^T a(X_s)ds = \mu. 
$$
Следовательно, выполнено равенство 
\begin{equation}\label{rmu}
r= \mu.
\end{equation}

~

{\bf 2 способ.} Пусть $a(\cdot)$ -- произвольная эргодическая стратегия. Известно \cite{Hasm, Hasm2}, что в этом случае существует инвариантная вероятностная мера $\mu^a$ и имеет место Закон больших чисел (ЗБЧ) в виде 
\begin{equation}\label{lln}
\frac{\displaystyle \int_0^T a(X_s)ds}{T} \to \mathsf E^{a,inv} a(X_0), \; t\to\infty, \quad \mathsf P^a -{\text{п.н.},}
\end{equation}
где $\mathsf E^{a,inv} a(X_0) \equiv \int a(y)\mu^a(dy)$ -- интеграл по инвариантной мере, соответствующей данной стратегии, и также 
\begin{equation}\label{Eax}
\frac{\displaystyle \mathsf E^a\int_0^T a(X_s)ds}{T} \to \mathsf E^{a,inv} a(X_0), \quad T\to\infty.
\end{equation}
Если допустить, что $\mathsf E^a a(X_0)< \mu$, то в силу (\ref{lln}) процесс окажется транзиентным и будет уходить к $+\infty$. Действительно, за большое время $T$ положительный прирост $\mu T$ c вероятностью, близкой к 1, намного превысит отрицательный $\int_0^T a(X_s)ds$ (с минусом). Случай же $\mathsf E^a a(X_0)>\mu$ при заданных ограничениях на модель невозможен. Действительно, в этом случае процесс должен был бы уходить на $-\infty$ в силу того же ЗБЧ (\ref{lln}), но это запрещено условием (\ref{a0}). Итак, эргодичность влечет за собой равенство $\mathsf E^{a,inv} a(X_0)=\mu$. Но тогда также имеем (\ref{rmu}) в силу (\ref{Eax}).

\end{proof}
\begin{rRemark}
Отметим, что, конечно, с математической точки зрения достаточно одного доказательства\footnote{Известен аргумент А.Н. Колмогорова об  обязательности нескольких  доказательств в исторических исследованиях и достаточности одного в математике, почему он и предпочел выбрать математику, а не историю.}. Тем не менее, эти два способа получения результата показывают различные возможности, которые могут оказаться полезны в более сложных случаях\footnote{Тут, все же, уместно вспомнить в какой-то степени противоположный совет Леонарда Эйлера, который говорил, что если к решению задачи намечается несколько дорог или тропинок, то крайне желательно пройти по ним по всем.}. Еще один, более стандартный способ доказательства приведен в следующем разделе. Он включен для того, чтобы продемонстрировать, что тут нет никакого противоречия с этим более стандартным подходом на основе уравнения Беллмана. Само уравнение будет решено не вполне строго, используя догадку. Однако, проверка решения очевидна, и на его основе будет еще раз повторен анализ оптимальных стратегий.  
\end{rRemark}

\begin{rEx}
Показательны следующие примеры. В самом общем случае  напомним, что $M>\mu$. Выберем произвольно $x_0\ge 0$ и стратегию $a(\cdot)$ вида 
$$
a(x) = \mu\,1(x\le x_0) - M\,1(x > x_0).
$$
В условиях задачи имеем экспоненциальную эргодичность (см., например, \cite{Ver2000}), единственную инвариантную меру $\mu^a$, причем последняя интегрирует любой полином. 

~

\noindent
\textbf{Случай $M=2\mu$}. Выбираем 
$$
a(x) = 2\mu\,1(x > x_0).
$$ 
Тогда снос имеет вид
$$
\mu - a(x) = \mu\,1(x\le x_0) - \mu\,1(x > x_0).
$$ 
Это симметричный случай, когда, ``очевидно'', что процесс в среднем половину времени проводит  слева от $x_0$, а половину справа. Стало быть, на этой стратегии достигается средний выигрыш
$$
r= 0*1/2 + 2\mu*1/2=\mu.
$$
Хотя для обоснования того, что $r=\mu$ следующее соображение уже лишнее, тем не менее, тут несложно найти стационарную плотность с особенностью производной в точке $x_0$: она имеет вид 
$$
p(x) = p(x_0)(\exp(\mu (x-x_0))\ 1(x\le x_0) + \exp(-\mu (x-x_0))\ 1(x>x_0)), 
$$
где из условия нормировки находится $p(x_0)=\mu/2$. 

~

\noindent
\textbf{Случай $M>2\mu$}. Хотя вариант выбора $a(x) = M\ 1(x>x_0)$ кажется весьма естественным, однако, чтобы свести задачу к предыдущей\footnote{Как известно, математики обожают это делать.}, тут можно снова выбрать 
$$
a(x) = 2\mu\,1(x > x_0),
$$ 
тогда вновь снос имеет вид
$$
\mu - a(x) = \mu\,1(x\le x_0) - \mu\,1(x > x_0), 
$$ 
и, стало быть, опять
$$
r= 0*1/2 + 2\mu*1/2=\mu.
$$

~
 
\noindent 
\textbf{Случай $M=C\mu$, $C>1$}. Оказывается, не обязательно делать снос симметричным. Выберем 
$$
a(x) = M\ 1(x>x_0).
$$
Найдем стационарную плотность $p(\cdot)$. Пусть ее значение в точке $x_0$ равно $p(x_0)$. Уравнение стационарности на полуоси $x>x_0$, дополненное условием непрерывной склейки 
$$
p(x_0-) = p(x_0+)=:p(x_0), 
$$
имеет вид
\[ 
\mathcal{L}^* p(x) = \frac12\ p''(x) - (b(x) p(x))' = 0, \quad x\neq x_0,
\] 
где $b(x)=-(C-1)\mu$ при $x>x_0=0$ и $b(x)=\mu$ при $x<x_0=0$. 

~

1) При $x>x_0$ имеем, 
\[
p'' + \gamma_1 p' = 0, \quad x > x_0, \quad \gamma_1 = (C-1)\mu>0.
\] 
Как несложно найти, и также несложно проверить, решением является функция 
\[
p(x)=p(x_0)e^{-\gamma_1 (x-x_0)}, \quad x>x_0.
\]

~

2) При $x<x_0$, обозначая $\gamma_2=\mu$, аналогично находим, 
\[
p(x)=p(x_0)e^{\gamma_2 (x-x_0)}, \quad x<x_0. 
\] 

~

Из условия нормировки 
\[ 1= p(x_0)\left( \int_0^{\infty} e^{-\gamma_1 x} dx + \int_0^\infty e^{-\gamma_2 x} dx \right)= p(x_0) \left( \frac{1}{\gamma_1} + \frac{1}{\gamma_2}, \right) 
\]
вычисляем,
\[
p(x_0) 
= \frac{\gamma_1 \gamma_2}{\gamma_1 + \gamma_2}.
\]

~

Значит, стационарные вероятности находиться справа и слева от $x_0$ равны, соответственно, 
\[
p_+ = \int_0^{\infty}  p_0 e^{-\gamma_1 x} dx = \frac{p_0} {\gamma_1} = \frac{\gamma_1 \gamma_2}{\gamma_1 + \gamma_2} \cdot \frac{1}{\gamma_1} = \frac{\gamma_2}{\gamma_1 + \gamma_2} = \frac{\mu}{(C-1)\mu+\mu} = \frac{1}{C},
\]
и
\[
p_- = \frac{\gamma_1}{\gamma_1 + \gamma_2}=\frac{(C-1)\mu}{(C-1)\mu+\mu}=\frac{C-1}{C}.
\] 
Стало быть, средний выигрыш при данной стратегии равен
\[
p_- \cdot 0 + p_+ \cdot C\mu = \mu,  
\] 
что и требовалось проверить. 

~

Отметим еще, что при стратегии $a(x) = c \mu \ 1(x>x_0)$ с условием $1<c\le M/\mu$,  выигрыш тоже оказывается равен $\mu$. Действительно, это вытекает буквально из той же выкладки  с заменой $C$ на $c$.

\end{rEx}

\medskip

Еще раз подчеркнем, что для нахождения $r$, как и для построения оптимальной стратегии, не было использовано уравнение Беллмана. Однако, это не означает, что в данной модели его использовать  невозможно. Как можно это сделать -- будет продемонстрировано в следующем разделе.

\section{Об  уравнении Беллмана}

Для поиска оптимального управления можно использовать эргодическое уравнение Беллмана (см.  \cite{BB, Borkar, Ver}):
\begin{equation}
\sup_{a \in [0, M]} \left[ (\mu - a)V'(x) + \frac12 \ V''(x) + a - r \right] = 0. \label{HJB}
\end{equation}
Здесь решением является пара $(V(\cdot),r)$. Известно  (см. \cite{BB, Borkar, Ver}), что при наложенных ограничениях решение единственно в классе пар, где $r$ -- постоянная, а функция $V$ растет не быстрее некоторого полинома и имеет две соболевские производные, локально интегрируемые в $L_p$ c любым $p>1$. Кроме того, компонента $r$ тут определяется единственным образом, а вспомогательная функция $V$ -- с точностью до аддитивной константы. Далее это будет продемонстрировано непосредственно для данного примера, а пока отметим (догадаемся), что решением является пара $V(x) = x$, $r=\mu$, что легко поддается проверке: $x'\equiv 1$, и $\mu - a + a - \mu = 0$. 

Если решение $(V,r)$ получено, то оптимальная стратегия может быть найдена в виде
\begin{equation}\label{aopt}
\bar a(x) \in \mathop{\text{Argsup}}\limits_a\left[ (\mu - a)V'(x) + \frac12 \ V''(x) + a - r \right] =
\mathop{\text{Argsup}}\limits_a\left[ a(1-V'(x))\right].
\end{equation}

~

Обоснование корректности такой постановки и использования уравнения \eqref{HJB} можно найти в работе \cite{BB}, где рассматривается общая задача эргодического управления одномерными диффузиями с условием ``near-monotonicity'' на функцию стоимости. В частности, в Теореме 3.4 указанной статьи доказано, что если существует стационарное марковское управление, при котором достигается минимум средней стоимости, то:

\begin{itemize}
    \item существует единственное решение уравнения Беллмана для эргодической задачи;
    \item найденное стационарное управление является оптимальным;
    \item решение \( V(x) \) можно представить как предел решений дисконтированных задач при \( \alpha \to 0 \).
\end{itemize}

В нашей постановке функция стоимости \( f(a,x) = a \) является непрерывной и ограниченной. Более того, при \( M > \mu \) и при подходящей стратегии процесс является эргодическим, что обеспечивает наличие инвариантной меры, аналогично условию ``существования стабильного марковского управления'' из статьи \cite{BB}.

Следовательно, в рамках модели можно использовать уравнение  (\ref{aopt}) для поиска оптимального управления. 
Согласно (\ref{aopt}), оптимальная стратегия должна иметь следующую структуру:
\[
a^*(x) =
\begin{cases}
0, & \text{если } V'(x) > 1, x<x_0=0,\\
M, & \text{если } V'(x) < 1, x>0,\\
\text{любое } a \in [0, M], & \text{если } V'(x) = 1, x=0.
\end{cases}
\]
Однако, как мы видели,  $V'(x)\equiv 1$. Следовательно, оптимальной будет любая стратегия $\bar a(\cdot)$ со значениями из отрезка $[0,M]$, лишь бы функция $\bar a(\cdot)$ была борелекской и были бы выполнены дополнительные условия эргодичности процесса. 

~

Покажем (немного нестрогий) алгоритм решения уравнения Беллмана, которое собираемся теперь использовать для анализа оптимальных стратегий. При этом без ограничения общности можем и будем считать, что $x_0=0$. Также, считаем, что значение $r=\mu$ уже установлено, см. предыдущий раздел. 

Итак, при $x<0$ имеем уравнение
\[
V''(x) + \mu V'(x) - r = 0 \quad \text{т.е.} \quad V''(x) + \mu V'(x) = r, \quad x<0.
\]
Вычисляя частное решение $V_{s}(x) = \frac{r}{\mu}x$ неоднородного ОДУ  и общее решение $V_{g}(x) = C_1 + C_2 e^{-\mu x}$  однородного ОДУ ($V'' + \mu V' =0$), имеем, 
\[
 V_1(x) = C_1 + C_2 e^{-\mu x} + \frac{r}{\mu} x.
\]

Аналогично при $x>0$ имеем уравнение
\[
V'' + (\mu - M) V' + M - r = 0. 
\]
Его решение имеет вид
\[
V_2(x) = \tilde{C}_1 + \tilde{C}_2 e^{(M - \mu)x} + \frac{r - M}{\mu - M} x, \quad x>0.
\]

Из условий гладкой склейки\footnote{Эти условия не обсуждаются. Если с их помощью удается найти классическое или соболевское решение, к которому можно применять формулу Ито после подстановки диффузии в функцию $V$, то проверка того, что вторая компонента $r$ является решением задачи (\ref{sup}), является стандартной процедурой.} 
\[
V_1(0-) = V_2(0+) \quad \text{и} \quad V'_1(0-) = V'_2(0+),
\]
находим
\[ 
V_1(0-) = C_1 + C_2 = V_2(0+) = \tilde{C}_1 + \tilde{C}_2,
\]
и
\[
V_1'(x) = -\mu C_2 e^{-\mu x} + \frac{r}{\mu}, \quad x<0.
\]
Далее,
\[
V_1'(0-) = -\mu C_2 + \frac{r}{\mu} = 1, 
\]
и поскольку $r=\mu$, то заключаем, что $C_2 = 0$.
  
Аналогично,
\[
V_2'(x) = \tilde{C}_2 e^{(M-\mu)x} (M-\mu) + \frac{r - M}{\mu - M},
\] 
откуда
\[
V_2'(0+) = \tilde{C}_2 (M-\mu) + \frac{r - M}{\mu - M} =1, 
\]
и поскольку $\frac{r - M}{\mu - M}=1$,  то заключаем, что $\tilde{C}_2=0$.

Итак, $C_1 = \tilde{C}_1 = :C$. 
Стало быть, 
\[
V_1(x) = C + \frac{r}{\mu} x = C + x,
\]
и аналогично
\[
V_2(x) = C + x.
\]
Итак, окончательно находим
$$
V(x) = C+x.
$$
Тот факт, что постоянную $C$ здесь определить не удается, не должен удивлять, поскольку в уравнение Беллмана входят лишь первые две производные $V$, но не сама эта функция, так что первая компонента решения определяется с точностью до произвольной аддитивной постоянной.

\medskip

Подставляя найденную {\em вспомогательную} функцию $V$ в уравнение (\ref{aopt}), 
видим, что оно предписывает выбирать любую (борелевскую) функцию $a(x) \in [0,M]$. Однако, не следует забывать, что стратегия обязана быть эргодической, а также что у нас есть условие (\ref{a0}), так что, строго говоря, этот выбор не совсем произвольный. Все это полностью согласуется с ответом из  раздела 3.

\medskip

Для завершения рассмотрения данной модели напомним, что при любой ограниченной борелевской стратегии $a(\cdot)$ уравнение (\ref{sde1}) имеет единственное сильное решение, в силу результатов любой из работ\footnote{Эти работы все, конечно, весьма разные, однако, есть общее пересечение, которое и позволяет в данном случае на любую из них сослаться.} \cite{Nakao, Zvon, Ver79, Ver80}.

\section*{Выводы}

Проведённое исследование позволяет сделать некоторые  принципиальные выводы о структуре оптимального управления и поведении системы в стационарном режиме.

Во-первых, оказалось не сложно, из весьма общих соображений и не прибегая к уравнению Беллмана, установить равенство $r=\mu$. Интуитивно, это равенство напрямую связано с устойчивостью системы: если бы \( r > \mu \), то в среднем из резерва  изымалось бы больше капитала, чем поступает за счет премий;  а значит, величина резерва в среднем убывала бы, противореча стационарности. Аналогично, \( r < \mu \) не соответствует оптимальности, так как всегда можно увеличить изъятие, не нарушая баланс и сохраняя эргодичность. Следовательно, из всех допустимых (то есть, таких, при которых существует сильное решение уравнения (\ref{sde1})) стратегий оптимальными являются те и только те, которые эргодичны в смысле условия (\ref{Xt0}) и обеспечивают равенство \( r = \mu \).

Во-вторых, 
построение стационарной плотности при кусочно-постоянных стратегиях управления учитывает, что при \( x < 0 \) или при \(x<x_0\) с некоторым $x_0\ge 0$, компания ничего не изымает (управление \( a(x) = 0 \)), а при \( x > 0 \) или при  \(x>x_0\) может изымать максимальную сумму \( a(x) = M >\mu\). При этом дрейф процесса слева от $x_0$ равен \( \mu \), а справа  $-(M-\mu)$, что обеспечивает возвращение процесса в окрестность точки \( x = 0 \) и наличие эргодичности с экспоненциальной сходимостью к стационарному режиму. Возможны и другие стратегии, достигающие той же цели, лишь бы процесс при них оказался эргодическим с выполнением условия (\ref{Xt0}). Любопытным можно назвать тот факт, что и выбор параметра $x_0\ge 0$ здесь произволен.

И, в третьих, уравнение Беллмана в данном примере допускает сравнительно несложное решение (учитывающее уже найденное равенство $r=\mu$), и его  анализ подтверждает выводы об оптимальных стратегиях, полученные без его использования.



\begin{thebibliography}{99}

\bibitem{Taksar2000}
M.I. Taksar, Optimal Risk and Dividend Distribution Control Models for an Insurance Company, Mathematical Methods of Operations Research, 2000, 51,  1--42. https://doi.org/10.1007/s001860050001

\bibitem{Ver2000}
А.Ю. Веретенников, О полиномиальном перемешивании и скорости сходимости для стохастических разностных и дифференциальных уравнений. Теория вероятн. и ее применен., 1999, 44(2), 317--321. http://mi.mathnet.ru/tvp766. 

\bibitem{Borkar}
A. Arapostathis, V.S. Borkar, M.K. Ghosh, 
Ergodic Control of Diffusion Processes, 
Cambridge University Press, 2011. 
https://doi.org/10.1017/CBO9781139003605

\bibitem{Ver}
S.V. Anulova, H. Mai, A.Yu. Veretennikov, On Iteration Improvement for Averaged Expected Cost Control for One-Dimensional Ergodic Diffusions, SIAM Journal on Control and Optimization. 2020, 58(4), 2312--2331. doi 10.1137/19M1271944 

\bibitem{Hasm}
Р.З. Хасьминский, Эргодические свойства возвратных диффузионных процессов и стабилизация решений задачи Коши для параболических уравнений, Теория вероятн. и ее примен., 1960, 5(2),  196--214. https://www.mathnet.ru/rus/tvp4825.
%

\bibitem{Hasm2}
Р.З. Хасьминский, 
Устойчивость систем дифференциальных уравнений при случайных возмущениях, М., Наука, 1969.


\bibitem{BB}
A. Bensoussan, V. Borkar, Ergodic control problem for one-dimensional diffusions with near-monotone cost, Systems \& Control Letters, 1984, 5(2), 127--133. https://doi.org/10.1016/0167-6911(84)90021-5


\bibitem{Nakao}
S. Nakao, On the pathwise uniqueness of solutions of one-dimensional stochastic differential equations, Osaka J. Math. 1972, 9, 513--518. 
https://projecteuclid.org/journals/osaka-journal-of-mathematics/volume-9/issue-3/On-the-pathwise-uniqueness-of-solutions-of-one-dimensional-stochastic/ojm/1200693939.pdf

\bibitem{Zvon}
А.К. Звонкин, Преобразование фазового пространства диффузионного процесса, уничтожающее снос, Матем. сб., 1974, 93(135):1, 129--149. 
https://www.mathnet.ru/rus/sm2963


\bibitem{Ver79}
А.Ю. Веретенников,
О сильных решениях стохастических дифференциальных уравнений, ТВП, 24:2 (1979), 348–360. Engl. transl. On the Strong Solutions of Stochastic Differential Equations, Theory Probab. Appl. Volume 24, Issue 2, pp. 354--366.
https://www.mathnet.ru/rus/tvp2867

\bibitem{Ver80}
А.Ю. Веретенников, 
О сильных решениях и явных формулах для решений стохастических интегральных уравнений, Матем. сб., 1980, 111(153):3, 434--452. 
https://doi.org/10.1137/1124039

\end{thebibliography}
\end{document}